\documentstyle[12pt]{amsart}

\makeatletter
\def\cases#1{\left\{\,\vcenter{\normalbaselines\m@th
    \ialign{$##\hfil$&\quad##\hfil\crcr#1\crcr}}\right.}
\makeatother

\def\R{{\bold R}}
\def\F{{\bold F}}

\author{Jie Wu}
\address{MSRI\\
1000 Centennial Drive\\
Berkeley, CA 94720-5070}
\title{On the Homology of Configuration Spaces $ C((M,M_o) \times \R^n; X)$}
\thanks{Research at MSRI is supported in part by NSF grant DMS-9022140}

\newtheorem{theorem}{Theorem}
\newtheorem{lemma}{Lemma}[section]
\newtheorem{proposition}[lemma]{Proposition}

\newtheorem{otherlemma}[lemma]{Lemma}\newtheorem{remark}[lemma]{Remark} 

\begin{document}

\maketitle
\begin{abstract}
The homology with coefficients in
a field  of the  configuration spaces $C(M\times \R ^n,M_o\times \R ^n;X)$ is determined in this paper.
\end{abstract}

\section{ Introduction.}

The purpose of this paper is to determine the homology (with coefficients in
a field $\F $) of configuration spaces $C((M,M_o)\times \R ^n;X)$.
The homology groups  
$$H_{*}(C((M,M_o)\times \R ^n;X);\F )$$ are
determined by dim $M$, $H_{*}(M,$ $M_o;$ $\F )$ and $H_{*}(X;\F )$.
This answers a conjecture of F. Cohen and L. Taylor [10].

Let $A$ be a topological space. The ordered configuration space $\widetilde{C%
}^k(A)$ is the subspace of $A^k$ consisting of all $k$-tuples $(a_1,...,a_k)$
such that $a_i\neq a_j$ for $i\neq j$, where $A^k=A\times ...\times A$ with $%
k$ copies. Define the configuration space 
$$
C(A,A_o;X)=\coprod_{k=1}^\infty \widetilde{C}^k(A)\times _{\Sigma _k}X^k 
/\approx 
$$
where $A_o$ is a closed subspace of $A$, $X$ is a space with
non-degenerate base-point $*$, and $\approx $ is generated by 
$$
(a_1,...,a_k;x_1,...,x_k)\approx (a_1,...,a_{k-1};x_1,...,x_{k-1}) 
$$
if $a_k\in A_o$ or $x_k=*$.

The spaces in the title are given by $C((A,A_o)\times \R ^n;X)=$ $C(A\times \R ^n,$ $A_o\times 
\R ^n;$ $X)$.
The applications of configuration spaces can be found in  [1, 3, 6, 7, 10, 11, 14, 17, 18]. The homology of configuration spaces for various cases can be found in [4, 5, 8, 9, 10, 11, 12, 13, 20]. We will consider the case that $A$ is a manifold and $A_{0}$ is a submanifold and $n>0$.

In this paper, a space $X$ will always mean a compactly generated weak
Hausdorff space with nondegenerate base-point $*$ such that $(X,*)$ is an
NDR-pair. Homology groups are always taken with coefficients in the field $\F $;$S_{*}(-)$ willalways mean the singular chain complex with coefficients in the field ${\bf F}$. A manifold will always mean smooth triangulated manifold. For each finitedimensional graded $\F $-module $V_{*}$ and  a connected space $X$, define
$$
{\cal C}^m(V_{*};X)=\bigotimes_{q=0}^m\bigotimes^{\beta
_q(V_{*})}H_{*}(\Omega ^{m-q}S^mX) 
$$
where $\beta _q(V_{*})=\dim {}_{\F }(V_{q})$ is the $q$-th Betti number of 
$V_{*}$ and $\Omega ^{m-q}S^mX=*$ if $q>m$.

If $V_q=0$ for all $q\geq m$, then each term $H_{*}(\Omega ^{m-q}S^mX)$ in $%
{\cal C}^m(V_{*};X)$ is an algebra with an algebra filtration induced by  $\Omega ^{m-q}S^mX\simeq C(\R ^{m-q},S^qX)$. This filters ${\cal C}^m(V_{*};X)$ by the tensor product filtration. The natural filtration of configuration spaces will be given in section 2. Now our main theorems are as follows

\begin{theorem}
Let $M$ be a smooth triangulated compact manifold and let $M_o$ be a smooth
triangulated compact submanifold of $M$,then

\begin{description}
\item[(1)]  for each simply connected space $X$ and $n\geq 1$, there is an $%
\F $-filtered {\bf module} isomorphism%
$$
\theta :{\cal C}^{\dim   M +n}(H_{*}(M,M_o);X)  %
\longrightarrow   H_{*}C((M,M_o)\times \R ^n,X) 
$$

\item[(2)]  for each simply connected space $X$ and $n\geq 2$, there is an $%
\F $-filtered {\bf algebra} isomorphism%
$$
\theta :{\cal C}^{\dim   M  +n}(H_{*}(M,M_o);X)  %
\longrightarrow   H_{*}C((M,M_o)\times \R ^n;X) 
$$
\end{description}
\end{theorem}
  
\begin{theorem}
Let $M$ be a smooth triangulated compact manifold and let $M_o$ be a smooth
triangulated compact submanifold of $M$, then there exist isomorphisms%
$$
\overline{H}*C((M,M_o)\times \R ^n;X)  \approx   %
\bigoplus_{k=1}^\infty \sigma ^{-2k}{\cal D}_k^{m+n}(H_{*}(M,M_o);S^2X) 
$$
as $\F $-modules for any space X and $n\geq 1$, where $m=$dim $M$, $%
{\cal D}_k^{m+n}(V_{*};X)=$%
$$
F_k{\cal C}^{m+n}(V_{*};X)/F_{k-1}{\cal C}^{m+n}(V_{*};X) 
$$
and $\sigma ^{-t}$ denotes the $t$-th desuspension of graded modules.
\end{theorem}

The article is organized as follows. In section 2, 
some basic properties of configuration spaces are studied. some product decompositions of certain configuration spaces are given in section 3. The proofs of Therem
A and Theorem B will be given in section 4. The author would like to thank  Professors Xueguang Zhou and Fred Cohen for their many helpful discussions and encouragement during the writing of the manuscript.The author also would like to thank Miguel Xicotencatl for his help to type the manuscript.

\bigskip
\bigskip

\section{Basic properties of configuration spaces.}

In this section,  some basic properties of $%
C(A,A_o;X)$ are recalled. Given an embedding $(A,A_o)$ $\rightarrow (A^{\prime
},A_o^{\prime })$ of space-pairs and a pointed map $X\rightarrow X^{\prime }$%
, there is an induced map $C(A,A_o;X)\longrightarrow C(A^{\prime
},A_o^{\prime };X^{\prime })$. Hence the homotopy type of $C(A,A_o;X)$ is an
invariant of the (relative) isotopy type of $(A,A_o)$ and the homotopy type
of $(X,*)$. The length of a configuration induces a natural filtration
of $C(A,A_o;X)$ by the closed subspaces%
$$
F_kC=F_kC(A,A_o;X)=\coprod_{j=1}^k\widetilde{C}^j(A)\times _{\Sigma _j}X^j 
  /\approx 
$$
$F_oC=*$, $F_1C=A/A_o\wedge X$, see [3]. It is easy to see that $C(A\amalg
A^{\prime },A_o\amalg A_o^{\prime };$ $X)$ is homeomorphic to the product
$$C(A,A_{o};X) \times C(A^{\prime},A_{o}^{\prime};X)$$ and the homeomorphism preserves
the filtration. Hence $C(A,A_o;X)$ is a filtered H-space if there is an
embedding $e:(A\amalg A,A_o\amalg A_o)\longrightarrow (A,A_o)$ such that $%
e|_{\makebox{first copy of }(A,A_o)}$ and $e|_{\makebox{second copy of }(A,A_o)}$
are (relatively) isotopic to the identity map of $(A,A_o)$. In particular, $%
C((A,A_o)\times \R ^n;X)$ is a filtered H-space for each $n\geq 1$.

Define
$$
\widetilde{C}^k(A|A_o)=\{(a_1,...,a_k)\in \widetilde{C}^k(A)  |\makebox{
some }a_j\in A_o\} 
$$

and 
$$
(X|*)^k=\{(x_1,...,x_k)\in X^{k  }|\makebox{ some }x_j=*\}\makebox{ .} 
$$

If $(A,A_o)$ is a relative CW-complex, there are $\Sigma _k$-equivalent
cofibrations $\widetilde{C}^k(A|A_o)$ $\rightarrow $ $\widetilde{C}^k(A)$
and $(X|*)^k\rightarrow X^k$ and therefore there are cofibrations%
$$
\widetilde{C}^k(A|A_o)\times _{\Sigma _k}X^k  \cup   \widetilde{C%
}^k(A)\times _{\Sigma _k}(X|*)^k  \longrightarrow   \widetilde{C}%
^k(A)\times _{\Sigma _k}X^k 
$$
and $F_{k-1}C\rightarrow F_k$ with the same cofibre $%
D_k(A,A_o;X)=F_kC/F_{k-1}C$ $\cong $ 
$$
\widetilde{C}^k(A)\times _{\Sigma _k}X^{k}/ 
(\widetilde{C}^k(A|A_o)\times _{\Sigma _k}X^k\cup \widetilde{C}^k(A)\times
_{\Sigma _k}(X|*)^k)
$$
see [2, pp 231-239; 14 pp 162-172 and Thm. 7.1]. We call $D_k(A,A_o;X)$ the $%
k$-adic construction of $C(A,A_o;X)$.

\begin{proposition}
If $(A,A_o)$ is a relative CW-complex and $k\geq 1$, then there is an
isomorphism of $\F $-modules%
$$
\overline{H}_{*}{\cal D}_k(A,A_o;X)\longrightarrow H_{*}\left( S_{*}( 
\widetilde{C}^k(A),\widetilde{C}^k(A|A_o))\otimes _{\Sigma _k}(\overline{H}%
_{*}(X))^{\otimes k}\right) 
$$
\end{proposition}

The following lemma is useful.

\begin{otherlemma}
Let $C$ be a free $\Sigma _k$-chain complex and let $K$ and $L$ be chain
complexes. If the chain maps $f,g:K\rightarrow L$ are homotopic, then $%
1\otimes f^{\otimes k}$ and $1\otimes g^{\otimes k}:C\otimes L^{\otimes
k}\longrightarrow C\otimes K^{\otimes k}$ are $\Sigma _k$-equivariently
homotopic, where $\Sigma _k$ acts diagonally on $C\otimes L^{\otimes k}$ and 
$C\otimes K^{\otimes k}$.
\end{otherlemma}

{\em Proof.} Let I be the unit chain complex with $0$-simplexes $
\overline{0}$ and $\overline{1}$ and $1$-simplex $\overline{I}$ and the
differential $\partial (\overline{I})=\overline{1}-\overline{0}$, and let $%
D:I\otimes K\rightarrow L$ the chain homotopy between $f$ and $g$, the
composite of $\Sigma _k$-equivariant chain maps%
$$
I\otimes C\otimes K^{\otimes k}\stackrel{\varphi \otimes 1}{\longrightarrow }%
I^{\otimes k}\otimes C\otimes K^{\otimes k}\approx C\otimes (I\otimes
K)^{\otimes k}\stackrel{1\otimes D^{\otimes k}}{\longrightarrow }C\otimes
L^{\otimes k} 
$$
is the required $\Sigma _k$-equivariant chain homotopy, where the $\Sigma _k$%
-equivariant chain map $\varphi :I\otimes C\rightarrow I^{\otimes k}\otimes
C $ is defined as follows.

Let ${\cal B}$ be a $\Sigma _k$-basis of $C$, define $\varphi (\overline{0}%
\otimes c)=\overline{0}^k\otimes c$, $\varphi (\overline{1}\otimes c)= 
\overline{1}^k\otimes c$ for each$c\in {\cal B}$ and define%
$$
\varphi (\overline{I}\otimes \sigma \cdot c)=\sigma \cdot (\overline{0}%
^{k-1}\otimes \overline{I}+\sum_{j=1}^{k-2}\overline{0}^{k-1-j}\otimes 
\overline{I}\otimes \overline{1}^j+\overline{I}\otimes \overline{1}%
^{k-1})\otimes \sigma \cdot c 
$$
for each  $\sigma \in \Sigma _k$.

\noindent{\em Proof of Proposition 2.1.}  By Eilenberg-Zilber Theorem,
there are isomorphisms 
$$
\overline{H}_{*}D^k(A,A_o;X)\approx H_{*}(S_{*}(\tilde C^k(A),\tilde
C^k(A|A_o))\otimes _{\Sigma _k}S_{*}(X^k,(X|*)^k)) 
$$
$$
\approx H_{*}(S_{*}(\tilde C^k(A),\tilde C^k(A|A_o))\otimes _{\Sigma
_k}(S_{*}(X,*))^{\otimes k}). 
$$
Since $S_{*}(\tilde C^k(A),\tilde C^k(A|A_o)$ is a free $\Sigma _k$-chain
complex over $\F $ and $S_{*}(X,*)\simeq \overline{H}_{*}(X)$ as chain
complexes over $\F $, where $\overline{H}_{*}(X)$ with trivial
differential (see[19, Lemma VIII.3.1]). The assertion follows by the above
Lemma.

\begin{proposition}
If $(A,A_o)$ is a relative CW-complex, then%
$$
\overline{H}_{*}D_k(A,A_o;X)\approx \sigma ^{-2k}\overline{H}%
_{*}D_k(A,A_o;S^2X) 
$$
where $\sigma ^{-t}$ is the $t$-th desuspension.
\end{proposition}

{\em Proof.} By Proposition 2.1, there are isomorphisms%
$$
\overline{H}_{*}D^k(A,A_o;X)\approx H_{*}(S_{*}(\tilde C^k(A),\tilde
C^k(A|A_o))\otimes _{\Sigma _k}(\overline{H}_{*}(X))^{\otimes k}) 
$$
and%
$$
\overline{H}_{*}D^k(A,A_o;S^2X)\approx H_{*}(S_{*}(\tilde C^k(A),\tilde
C^k(A|A_o))\otimes _{\Sigma _k}(\overline{H}_{*}(S^2X))^{\otimes k}) 
$$
$$
\approx \sigma ^{2k}H_{*}(S_{*}(\tilde C^k(A),\tilde C^k(A|A_o)\otimes
_{\Sigma _k}\overline{H}_{*}(X)^{\otimes k}) 
$$

\begin{proposition}
If $(A,A_o)$ is a relative CW-complex with an embedding 
$$\coprod_{k=1}^%
\infty \tilde C^k(A)/\Sigma _k\rightarrow R^\infty ,$$ then there is a stable
equivalence%
$$
\sigma :C(A,A_o;X)\longrightarrow \bigvee_{k=1}^{\infty} D_k(A,A_o;X) 
$$
via stable equivalences%
$$
\sigma _k:F_kC(A,A_o;X)\longrightarrow \bigvee_{j=1}^k D_j(A,A_o;X) 
$$
\end{proposition}

{\em Proof.}  If $(A,A_o)=(M,M_o)$ a manifold-pair, this is
proved in [3, Prop. 3]. For the general cases, one proceeds in the same way.

\begin{proposition}
Let $X$ be a path connected space and let $(A,A_o)$ be a relative CW-complex with
an embedding $\coprod_{k=1}^\infty \tilde C^k(A)/\Sigma _k\longrightarrow
R^\infty $. There exists an isomorphism of $\F $-modules%
$$
\varphi :H_{*}C(A,A_o;\bigvee_{\alpha \in I}S^{n_\alpha })\longrightarrow
H_{*}C(A,A_o;X) 
$$
where $\{n_\alpha |\alpha \in I\}$ is determinined by $\overline{H}_{*}(X)$.
\end{proposition}

{\em Proof.} Let $\{x_\alpha |\alpha \in I\}$ be a basis of $
\overline{H}_{*}(X)$ and let $n_\alpha =|x_\alpha |$ for $\alpha \in I$, we
have an $\F $-isomorphism%
$$
\varphi _1:\overline{H}_{*}(\vee_\alpha S^{n_\alpha })\rightarrow 
\overline{H}_{*}(X) 
$$
Now the assertion follows by Proposition 2.1 and 2.3.

\begin{proposition}
Let $M$ be a smooth compact manifold and let $M_o$ and $N$ be the smooth
compact submanifolds of $M$ with codim$N=0$. If $N/M_o\cap N$ or $X$ is path
connected, then 
$$
C(N,N\cap M_o;X)\rightarrow C(M,M_o;X)\rightarrow C(M,N\cup M_o;X) 
$$
is a quasifibration. 
\end{proposition}

[5, pp. 113].

\bigskip

Now let $M$ be an $m$-manifold and let $W$ be a $m$-manifold without
boundary which contains $M$, e.g. $W=M$ if $M$ is closed, or $W=M\cup
\partial M\times [0,1)$ if $M$ has boundary. Let $\xi $ be the principal $%
O(m)$-bundle of the tangent bundle of $W$. Let $\Gamma _{\xi
[S^mX]}(B,B_o)$ be the space of cross sections of $\xi [S^mX]$ which are
defined on $B$ and take values at $\infty \wedge X$  on $B_o$ for each subspace pair $(B,B_o)$ in $W$%
, where $\xi [S^mX]$ is the associated bundle and $O(m)$ acts diagonally on $%
S^mX=S^m\wedge X$, trivially on $X$ and canonically on $S^m\cong R^m\cup
\{\infty \}$.

\begin{proposition}
Let $M$ be a smooth compact manifold and let $M_o$ be a smooth compact
submanifold of $M$. If $M/M_o$ or $X$ is path connected, then there is a (weak)
homotopy equivalence%
$$
C(M,M_o;X)\rightarrow \Gamma _{\xi [S^mX]}(W-M_o,W-M) 
$$

\end{proposition}

[5, Proposition 3.1 and 3, pp. 178]

\bigskip

\begin{remark}
By Proposition 2.6,there is a homotopy equivalence 
$$
C((M,M_o)\times \R ^n;X)\simeq \Omega ^nC(M,M_o;S^nX) 
$$
if $M/M_o$ or $X$ is path connected.
\end{remark}
\section{Decomposition Theorems}

In this section, assume that $(M,M_o)$ is a smooth triangulated
compact manifold-pair with $m=$dim$M$ and $W$ is a smooth $m$-manifold
without boundary which contains $M$ and $\xi _W$ is the principal $O(m)$%
-bundle of the tangent bundle of $W$ (see Proposition 2.6). Let $\overline{W}%
=M$ if $M$ is closed, or $\overline{W}=M\times [0,\frac 12]$ if $M$ has
boundary.

\begin{lemma}
If $X$ is path connected, then there is a (weak) homotopy equivalence%
$$
C((M,M_o)\times \R ;X)\rightarrow \Gamma _{\xi _W[\Omega
S^{m+1}X]}(W-M_o,W-M) 
$$
where $O(m)$ acts on $\Omega S^{m+1}X$ via the functor $\Omega S(-)$,i.e. by taking the homeomorphism $\Omega S(\sigma \wedge id|_{X}): \Omega S^{m+1}X \longrightarrow \Omega S^{m+1}X$ for each $\sigma \in O(m)$.
\end{lemma}

{\em Proof.} Notice that $(M,M_o)\times \R $ is isotopic to $%
(M,M_o)\times I$. By Proposition 2.6, it is easy to see that%
$$
C((M,M_o)\times \R ;X)\simeq \Gamma _{\xi _{W\times \R %
}[S^{m+1}X]}((W-M_o,W-M)\times (I,\partial I)) 
$$
Since $\xi _{W\times \R }[S^{m+1}X]=\pi ^{*}\xi _W[S^{m+1}X]$, where $%
\pi :W\times \R \rightarrow W$ is the projection, there is%
$$
\Gamma _{\xi _{W\times \R }[S^{m+1}X]}((W-M_o,W-M)\times (I,\partial I)) 
$$
$$
\simeq \Gamma _{\xi _W[\Omega S^{m+1}X]}(W-M_o,W-M) 
$$

\begin{lemma}
Let $(N,N_o)$ be a compact submanifold-pair in $\overline{W}$ such that $%
N\subseteq \overline{W}-\partial \overline{W}$, then there is a (weak)
homotopy equivalence%
$$
C((\overline{W}-\nu (N_o),\overline{W}-\nu (N))\times \R ;X)\rightarrow
\Gamma _{\xi _W[\Omega S^{m+1}X]}(N,N_o) 
$$
if $X$ is path connected, where$(\nu (N),\nu (N_o))$ is an open tubular
neighborhood of $(N,N_o)$ in $\overline{W}-\partial \overline{W}$.
\end{lemma}

{\em Proof.} By Lemma 3.1, there is a homotopy equivalence%
$$
C((\overline{W}-\nu (N_o),\overline{W}-\nu (N))\times \R ;X)\simeq
\Gamma _{\xi _W[\Omega S^{m+1}X]}(W-(\overline{W}-\nu (N)),W-(\overline{W}%
-\nu (N_o))) 
$$
$$
\cong \Gamma _{\xi _W[\Omega S^{m+1}X]}(\nu (N),\nu (N_o)) 
$$
$$
\simeq \Gamma _{\xi _W[\Omega S^{m+1}X]}(N,N_o) 
$$

\begin{lemma}
Let $(N,N_o)$ be a submanifold-pair in $W$. There is a homeomorphism%
$$
\Gamma _{\xi _W[\Omega S^{mn+1}X]}(N,N_o)\cong \Gamma _{\xi _{W^n}[\Omega
S^{mn+1}X]}(\Delta _n(N),\Delta _n(N_o)) 
$$
for $n\geq 1$, where $\Delta _n:W\rightarrow W^n=W\times ...\times W$ the
diagonal map, $O(mn)$ acts on $\Omega S^{mn+1}X$ via the functor $\Omega S(-)
$ ,see Lemma 3.1 as above, and $O(m)$ acts on $\Omega S^{mn+1}X$ via the diagonal inclusion $%
O(m)\rightarrow O(mn),\sigma \rightarrow diag(\sigma,\sigma,...,\sigma)$.
\end{lemma}

{\em Proof.} The diagonal map $\Delta _n:W\rightarrow W^n$ induces a bundle map $
\overline{\Delta }_n:\xi _W\rightarrow \xi _{W^n}$. Consider the induced map on the total spaces of bundles 
$\overline{\Delta }_n\times 1:E\xi _W\times \Omega S^{mn+1}X\longrightarrow
E\xi _{W^n}\times \Omega S^{mn+1}X$.

For each (x,y)$\in E\xi _W\times \Omega S^{mn+1}X$ and $\sigma \in O(m)$, notice that  
$$
\overline{\Delta }_n\times 1(x\sigma ,\sigma ^{-1}y)=((x,\ldots ,x)\tilde
\sigma ,\Omega S(\sigma^{-1}\wedge \ldots \wedge \sigma^{-1}\wedge 1_{X})y) 
$$
$$
=((x,\ldots ,x)\tilde \sigma ,\tilde \sigma ^{-1}y) 
$$
where $\tilde \sigma =\left( 
\begin{array}{ccc}
\sigma &  &  \\  
& \ddots &  \\  
&  & \sigma 
\end{array}
\right) \in O(mn)$.

Hence $\overline{\Delta }_n\times 1$ induces a map%

$$
\Delta : E(\xi _W[\Omega S^{mn+1}X])=E\xi _W\times _{O(m)}\Omega S^{mn+1}X  \longrightarrow  E\xi
_{W^n}\times _{O(mn)}\Omega S^{mn+1}X   
$$
$$
 =E(\xi _{W^n}[\Omega S^{mn+1}X]). 
$$
Furthermore $\Delta$ induces a homeomorphism%
$$
\Gamma _{\xi _W[\Omega S^{mn+1}X]}(N,N_o)\cong \Gamma _{\xi _{W^n}[\Omega
S^{mn+1}X]}(\Delta _n(N),\Delta _n(N_o)), 
$$
since the induced map is 1-1,onto and open.%

The following lemma follows from the naturality of the Samelson products. 

\begin{lemma}
There exists $O(m)$-map $\Phi _n: (S^{m}X)^{(n)}\longrightarrow \Omega
S^{m+1}X$ which represents the Samelson products $[[E,E],\ldots ,E]$ for each $n\geq 2$, where $X^{(n)}$ is the reduced join of $n$ copies of $X$ and $E: S^{m}X \longrightarrow \Omega S^{m+1}X$ is the suspension.
\end{lemma}

Now we give some decomposition theorems.

\begin{theorem}
Let $(M,M_o)$ be a smooth triangulated compact manifold-pair, $n+m+1$ even
and $n\geq 1$. There exist a manifold-pair $(\tilde M,\tilde M_o)$ so that $%
\dim \tilde M=2m$ and 
$$
C((M,M_o)\times \R ;S^n)_{(p)}\stackrel{w}{\simeq }C((\tilde M,\tilde
M_o)\times \R ,S^{2n})_{(p)}\times C(M,M_o;S^n)_{(p)} 
$$
where $m=\dim M$ and $p$ is an odd prime or zero and $\stackrel{w}{\simeq}$ means the (weak) homotopy equivalence
\end{theorem}

{\em Proof.} By Lemma 3.1, there is a homotopy equivalence 
$$
C((M,M_o)\times \R ;S^n)\simeq \Gamma _{\xi _W[\Omega
S^{m+1+n}]}(W-M_o,W-M) 
$$

The inclusion $j:S^{m+n}\rightarrow \Omega S^{m+1+n}$ is an $O(m)$-map. By
Lemma 3.4, there is an $O(m)$-Samelson product $[j,j]:S^{2m+2n}\rightarrow
\Omega S^{m+1+n}$, which induces an $O(m)$-H-map $\overline{[j,j]}:\Omega
S^{2m+1+2n}\rightarrow \Omega S^{m+1+n}$.

Now consider the $O(m)$-map%
$$
(\overline{[j,j]},j):\Omega S^{2m+1+2n}\times S^{m+n}\longrightarrow \Omega
S^{m+1+n}, 
$$
which induces a bundle map%
$$
\xi _W[\Omega S^{2m+1+2n}\times S^{m+n}]\longrightarrow \xi _W[\Omega
S^{m+1+n}] 
$$
and therefore a map%
$$
\Gamma _{\xi _W[\Omega S^{2m+1+2n}\times S^{m+n}]}(W-M_o,W-M)\longrightarrow
\Gamma _{\xi _W[\Omega S^{m+1+n}]}(W-M_o,W-M) 
$$
Notice that%
$$
\begin{array}{c}
\Gamma _{\xi _W[\Omega S^{2m+1+2n}\times S^{m+n}]}(W-M_o,W-M)\simeq \\ 
\simeq \Gamma _{\xi _W[\Omega S^{2m+1+2n}]}(W-M_o,W-M)\times \Gamma _{\xi
_W[S^{m+n}]}(W-M_o,W-M) 
\end{array}
$$
By Lemma 3.3%
$$
\begin{array}{c}
\Gamma _{\xi _W[\Omega S^{2m+1+2n}]}(W-M_o,W-M)\simeq \Gamma _{\xi
_{W^2}[\Omega S^{2m+1+2n}]}(\Delta (W-M_o),\Delta (W-M)) \\ 
\simeq C((\tilde M,\tilde M_o)\times \R ;S^{2n})\makebox{ (Lemma 3.2)} 
\end{array}
$$
where $\dim \tilde M=2m$. Since $(\overline{[j,j]},j)$ is homotpy equivalent after $p$-localization, the assertion follows by induction on the handle
decomposition of $M$.

\begin{theorem}
Let $(M,M_o)$ be a smooth triangulated compact manifold-pair. If
$X_1,\allowbreak \ldots
,\allowbreak X_k$ are connected spaces, there exist manifold pairs
$$(M_\omega ,M_{o\omega
})$$ such that $C((M,M_o)\times \R ;X_1\vee \ldots \vee X_k)$ is (weak)
homotopy equivalent to the (weak) product%
$$
\prod_\omega C((M_\omega ,M_{o\omega })\times \R ;\omega (X_1,\ldots
,X_k)) 
$$
where $\omega $ runs over all addmissible words in $x_1,\ldots ,x_k$ (the
addmissible words see [21, pp. 511-514]) and $\omega (X_1,\ldots
,X_k)=X_1^{(a_1)}\wedge \ldots \wedge X_k^{(a_k)}$ and $X^{(n)}$ is the
reduced join of $n$ copies of $X$ and $a_i$ is the number of occurences of $%
x_i$ in the word $\omega $.
\end{theorem}

{\em Proof.} Denote $X=X_1\vee \ldots \vee X_k$. By Lemma 3.1, there is a homotopy equivalence%
$$
C((M,M_o)\times \R ;X)\simeq \Gamma _{\xi _W[\Omega
S^{m+1}X]}(W-M_o,W-M) 
$$
Consider the inclusions 
$$
i_t:S^mX_t\rightarrow \Omega S^{m+1}X\makebox{\qquad for }1\leq t\leq k 
$$
By Lemma 3.4, there are $O(m)$-Samelson products 
$$
\omega (i_1,\ldots ,i_k):S^{l(\omega )\cdot m}\omega (X_1,\ldots
,X_k)\longrightarrow \Omega S^{m+1}X 
$$
for each word $\omega $, where $l(\omega )$, is the length of $\omega 
$, which induces an $O(m)$-H-map%
$$
\tilde \omega (i_1,\ldots ,i_k):\Omega S^{l(\omega )\cdot m}\omega
(X_1,\ldots ,X_k)\longrightarrow \Omega S^{m+1}X 
$$
Consider the $O(m)$-map%
$$
f=\prod_\omega \tilde \omega (i_1,\ldots ,i_k):\prod_\omega \Omega S^{{l%
}(\omega )\cdot m}\omega (X_1,\ldots ,X_k)\longrightarrow \Omega S^{m+1}X 
$$
which induces a bundle map%
$$
\xi _W[\prod_\omega \Omega S^{l(\omega )\cdot m}\omega (X_1,\ldots
,X_k)]\longrightarrow \xi _W[\Omega S^{m+1}X] 
$$
and therefore a map%
$$
\tilde f:\Gamma _{\xi _W[\prod_\omega \Omega S^{l(\omega )\cdot
m}\omega (X_1,\ldots ,X_k)]}(W-M_o,W-M) 
$$
$$
\longrightarrow \Gamma _{\xi W[\Omega S^{m+1}X]}(W-M_o,W-M). 
$$
By Hilton-Milnor Theorem, $f$ is a (weak) homotopy equivalence. By induction on the handle decomposition of $M$,$\tilde f$ is a (weak) homotopy quivalence.

Notice that%
$$
\Gamma _{\xi _W[\prod_\omega \Omega S^{l(\omega )\cdot m}\omega
(X_1,\ldots ,X_k)]}(W-M_o,W-M) 
$$
$$
\simeq \Gamma _{\xi _W[\prod_\omega \Omega S^{l(\omega )+1}\omega
(X_1,\ldots ,X_k)]}(\overline{W}-v(M_o),\overline{W}-v(M)) 
$$
$$
\cong \prod_\omega \Gamma _{\xi _W[\Omega S^{l(\omega )+1}\omega
(X_1,\ldots ,X_k)]}(\overline{W}-v(M_o),\overline{W}-v(M)) 
$$
$$
\cong \prod_\omega \Gamma _{\xi _{W^{l(\omega )}}[\Omega S^{l%
(\omega )+1}\omega (X_1,\ldots ,X_k)]}(\Delta _{l(\omega )}(\overline{%
W}-v(M_o)),\Delta _{l(\omega )}(\overline{W}-v(M))) 
$$
$$
\simeq \prod_\omega C((M_\omega ,M_{o\omega })\times \R ;\omega
(X_1,\ldots ,X_k))\makebox{\qquad (Lemma 3.2)} 
$$
where $\dim M_\omega =\dim W^{l(\omega )}=l(\omega )\cdot m$.

The assertion follows.

\section{Proofs of Theorems A and B.}

\begin{lemma}
Let $(N,N_o)\subseteq (M,M_o)$ be the smooth compact submanifold pair of $%
(M,M_o)$ with $\dim N=\dim M$ and let $X$ be a simply connected space. If $%
H_{*}(N,N_o)\rightarrow H_{*}(M,M_o)$ is onto, then 
$$
H_{*}C((N,N_o)\times \R ^n;X)\longrightarrow H_{*}C((M,M_o)\times \R %
^n;X) 
$$
is onto for $n\geq 1$.
\end{lemma}

{\em Proof.} Since $(M,M_o)\times \R ^n$ is isotopic to $%
(M\times I^{n-1},M_o\times I^{n-1})\times \R $, it is sufficient to show
that%
$$
H_{*}C((N,N_o)\times \R ;X)\longrightarrow H_{*}C((M,M_o)\times \R %
;X) 
$$
is onto.

\begin{itemize}
\item[Step1:]  Assume that $X=S^n$ with $n>1$.

If $m+1+n$ is odd or $\F $ is of characteristic $2$, this was proved in
[5, Th. A], where $m=\dim M$.

Now assume that $m+1+n$ is even and $\F $ is of characteristic $\not =2$%
. By Theorem C, there is a $p$-homotopy equivalence 
$$
C((M,M_o)\times \R ;S^n)\simeq C((\tilde M,\tilde M_o)\times \R %
;S^{2n})\times C(M,M_o;S^n) 
$$
and%
$$
C((N,N_o)\times \R ;S^n)\simeq C((\tilde N,\tilde N_o)\times \R %
;S^{2n})\times C(N,N_o;S^n) 
$$

where $\dim \tilde M+1+2n=\dim \tilde N+1+2n=2m+2n+1$

and $\dim M+n=\dim N+n=m+n$.

The assertion follows by [5, Th. A].

\item[Step 2:]  Assume that $X=S^{d_1}\vee \ldots \vee S^{d_k}$ with $d_j>1$.

By Theorem D, there is a (weak) homotopy equivalence%
$$
C((M,M_o)\times \R ;X)\simeq \prod_\omega C((M_\omega ,M_{o\omega
})\times \R ;\omega (X)) 
$$
where $\omega (X)=S^{d_1a_1}\wedge \ldots \wedge
S^{d_ka_k}=S^{d_1a_1+\cdots +d_ka_k}$.

The assertion follows from Step 1.

\item[Step 3:]  Assume $H_{*}(X)$ is a finite dimensional $\F $-vector
space. The assertion follows from Step 2 and Proposition 2.4.

\item[Step 4:]  General case.

There exist $X_\alpha $ such that $H_{*}(X)= 
\lim_\alpha H_{*}(X_\alpha )$ and $H_{*}(X_\alpha )$ is finite
dimensional $\F $-vector space for each $\alpha $. The assertion follows
from Step 3.
\end{itemize}

{\em Proof of Theorem A.} First we will prove the absolute
case $M_o=\phi $ by induction on on a handle decomposition of $M$. If M
is a disjoint union $M_1\amalg M_2$, then 
$$
C(M\times \R ^n;X)\cong C(M_1\times \R ^n;X)\times C(M_2\times {\bf R%
}^n;X) 
$$
Hence we can restrict to connected manifolds and start with $M=I^m$, the
assertion is obvious for $M=I^m$. Assume that the assertion holds for $M$
and $\overline{M}=M\cup D$ with $D\cong I^m$ a handle of index $q$ i.e. $%
D\cap M\cong I^{m-q}\times \partial I^q$, $q\geq 1$ since $M$ is connected.
There is a cofibration%
$$
M\stackrel{i}{\rightarrow }\overline{M}\stackrel{j}{\rightarrow }(\overline{M%
},M)\simeq (I^m,I^{m-q}\times \partial I^q)\simeq (S^q,*) 
$$
and alternative

I. $H_q(\overline{M})\stackrel{j_{*}}{\longrightarrow }H_q(S^q)$ is onto.

II. $H_q(\overline{M})\longrightarrow H_q(S^q)$ is zero.

\begin{itemize}
\item[Case I.]  Consider the quasifibration%
$$
C(M\times \R ^n;X)\longrightarrow C(\overline{M}\times \R ^n;X)%
\stackrel{C(j)}{\longrightarrow }C((\overline{M},M)\times \R ^n;X) 
$$
$$
\simeq \Omega ^{m+n-q}S^{m+n}X 
$$

Since $j_{*}$ is onto, $C(j)_{*}$ is onto and the Serre spectral
sequence for the quasifibration above collapses. Hence there is a short exact
sequence of Hopf algebras%
$$
H_{*}C(M\times \R ^n;X)\succ \longrightarrow H_{*}C(\overline{M}\times 
\R ^n;X)\longrightarrow \succ H_{*}\Omega ^{m+n-q}S^{m+n}X 
$$

By Proposition 2.3, there is a commutative diagram%
$$
\begin{array}{cccccccccccc}
\Sigma^{\infty}C(\overline{M}\times \R ^n;X)
&\stackrel{p}{\rightarrow}&\Sigma^{\infty}D_k\overline{M}\times\R ^n; X)\\
\downarrow  &   &\downarrow \\ 
\Sigma^{\infty}C((\overline{M}, M)\times \R ^n;X) &\stackrel{p}{\rightarrow}&\Sigma^{\infty}D_k((\overline{M},M)\times\R ^n; X) 
\end{array}
$$
for each $k\geq1$, where $p$ is the projection.

Thus\quad 
$$
\overline{H}_{*}D_k(\overline{M}\times \R ^n)\longrightarrow \overline{H}%
_{*}D_k((\overline{M},M)\times \R ^n;X) 
$$
is an epimorphism and $F_rH_{*}C(\overline{M}\times \R ^n;X)\longrightarrow
F_rH_{*}\Omega ^{m+n-q}S^{m+n}X$ is an epimorphism by Proposition 2.3.
\end{itemize}

Hence there exists an $\F $-map%
$$
\varphi :H_{*}\Omega ^{m+n-q}S^{m+n}X\longrightarrow H_{*}C(\overline{M}%
\times \R ^n;X) 
$$
such that

\begin{enumerate}
\item  $\varphi $ preserves the filtration.

\item  $C(j)_{*}\circ \varphi =id$
\end{enumerate}

Now the composite%
$$
H_{*}C(M\times \R ^n;X)\otimes H_{*}\Omega
^{m+n-q}S^{m+n}X\longrightarrow H_{*}C(\overline{M}\times \R %
^n;X)^{\otimes 2} 
$$
$$
\longrightarrow H_{*}C(\overline{M}\times \R ^n;X) 
$$
is an isomorphism of filtered modules and (1) follows.

If $n>1$, let $\varphi :QH_{*}(\Omega ^{m+n-q}S^{m+n}X)\longrightarrow
QH_{*}C(\overline{M}\times \R ^n;X)$ so that $QC(j)_{*}\circ \varphi =id$
and $\varphi $ preserves the filtration.

The homomorphism%
$$
QH_{*}\Omega ^{m+n-q}S^{m+n}X\longrightarrow QH_{*}C(\overline{M}\times {\bf %
R}^n;X)\longrightarrow H_{*}C(\overline{M}\times \R ^n;X)\  
$$
induces an algebra map%
$$
H_{*}\Omega ^{m+n-q}S^{m+n}X\longrightarrow H_{*}C(\overline{M}\times \R %
^n;X) 
$$
since $H_{*}C(\overline{M}\times \R ^n;X)$ is a commutative algebra and $%
H_{*}\Omega ^{m+n-q}S^{m+n}X$ is a free commutative algebra. Thus (2)
follows.

\begin{itemize}
\item[Case II.]  Consider the quasifibration%
$$
\Omega C((\overline{M},M)\times \R ^n;X)\longrightarrow C(M\times \R %
^n;X)\stackrel{C(i)}{\longrightarrow }C(\overline{M}\times \R ^n;X) 
$$
\end{itemize}

Since $i_{*}$ is onto, $C(i)_{*}$ is onto and the Serre spectral
sequence for the quasifibration above collapses. Hence%
$$
H_{*}C(\overline{M}\times \R ^n;X)\approx H_{*}C(M\times \R %
^n;X)//H_{*}\Omega ^{m+n-q+1}S^{m+n}X 
$$

Both (1) and (2) follow.

To treat the relative case, we can assume that $M_o$ is part of a closed
collar (see [5]). We will prove by induction om a handle decomposition of $M$%
, and start with $M=M_o$, the assertion is obvious for $M=M_o$ since $%
F(Mo,M_o)\simeq *$. Assume that the assertion holds for $(M,M_o)$. If $
\overline{M}=M\cup D$ with $D\cong I^m$ a hand of index $q$ in $M$. Clearly
we can assume that $q\geq 1$. There is a cofibration%
$$
(M,M_o)\rightarrow (\overline{M},M_o)\rightarrow (\overline{M},M)\simeq
(I^m,I^{m-q}\times \partial I^q) 
$$
and again an alternative:

III. $H_q(\overline{M},M_o)\rightarrow H_q(\overline{M},M)$ is onto.

IV. $H_q(\overline{M},M_o)\rightarrow H_q(\overline{M},M)$ is zero.

For Case III and IV, the assertion follows from Lemma 4.1 similar to Case I and II respectively.

{\em Proof of  Theorem B.}  By Theorem A, there is an isomorphism of filtered $\F $-modules%
$$
H_{*}C((M,M_o)\times \R ^n;S^2X)\approx {\cal C}%
^{m+n}(H_{*}(M,M_o);S^2X). 
$$

Hence%
$$
\overline{H}_{*}D_k((M,M_o)\times \R ^n;S^2X)\approx {\cal D}%
_k^{m+n}(H_{*}(M,M_o);S^2X)\makebox{.} 
$$

By Proposition 2.2 and Proposition 2.3, there are isomorphisms of $\F $-modules%
$$
\overline{H}_{*}D_k((M,M_o)\times \R ^n;X)\approx \sigma ^{-2k}\overline{%
H}_{*}D_k((M,M_o)\times \R ^n;S^2X) 
$$
and%
$$
\overline{H}_{*}C((M,M_o)\times \R ^n;X)\approx \bigoplus_{k=1}^\infty 
\overline{H}_{*}D_k((M,M_o)\times \R ^n;X) 
$$

The assertion follows.

\end{document}